\documentclass{article}
\usepackage{latexsym}

\title{Knot Theory of Coxeter type B and its physical
applications}
 \author{Reinhard H\"aring-Oldenburg\\
      Mathematisches Institut,
            Bunsenstr. 3-5, D-37073 G\"ottingen\\
 email: haering@cfgauss.uni-math.gwdg.de}

\newcommand{\RR}{\mathop{\rm I\! R}\nolimits}
\newcommand{\ZZ}{\mathop{\sl Z\!\!Z}\nolimits}

\newcommand{\deu}[1]{}
\newcommand{\eng}[1]{#1}

\newcommand{\kat}{{\cal C}}

\newcommand{\oalg}{{\cal A}}

\newcommand{\kegel}{{\cal O}}

\newcommand{\Mor}{{{\rm Mor}}}

\newcommand{\ev}{{{\rm ev}}}

\newcommand{\id}{{{\rm id}}}	  \newcommand{\tr}{{{\rm tr}}}

\newcommand{\Vec}{{{\rm Vec}}}
\newcommand{\Rep}{{{\rm Rep}}}

\newcommand{\BB}{{\rm B^\ast B}}
\newcommand{\ZA}{{{\rm ZA}}}
\newcommand{\ZB}{{{\rm ZB}}}

\newcommand{\HB}{{{\rm HB}}}

\newcommand{\dottedline}[5]{
                   \multiput(#1,#2)(#3,#4){#5}{\circle*{0.5}}
                        }

\newtheorem{de}{Definition}
\newtheorem{satz}{Proposition}

\begin{document}
\maketitle
\begin{abstract} Braid groups may be defined for every Coxeter diagram.
Artin's braid group is of type A. Analogs of Temperley-Lieb, Hecke and
Birman-Wenzl algebras exist for B-type. 

Our general hypothethis is that the braid group of B-type replaces  
Artin's braid group in most physical applications if the 
model is equipped with  a nontrivial boundary.
Solutions of a Potts model with a boundary and the reflection equation
illustrate this principle.

Braided tensor categories of  B-type and dually Coxeter-B braided
Hopf algebras are introduced. The occurrence of such categories in QFT on 
a half plane is discussed.
\end{abstract}

\sloppy

 \section{Introduction}
To every Coxeter diagram a braid group is associated that has the same 
presentation as the Coxeter group but has no  degree 2 relations
for the generators.  Artin's braid group is the braid group 
$\ZA_n$ of Coxeter type A. 
\begin{de} The braid group $\ZB_n$ of Coxeter type B is generated
by $X_0,X_1,\ldots,X_{n-1}$ with relations
\begin{eqnarray}
X_iX_j&=&X_jX_i\quad\mbox{if}\quad |i-j|>1\label{def1}\\
X_iX_jX_i&=&X_iX_jX_i\quad\mbox{if}\quad i,j\geq1,|i-j|=1\label{def2}\\
X_0X_1X_0X_1&=&X_1X_0X_1X_0\label{def3}
\end{eqnarray}
\end{de}
Generators $X_i,i\geq1$  satisfy the relations of Artin's braid group.

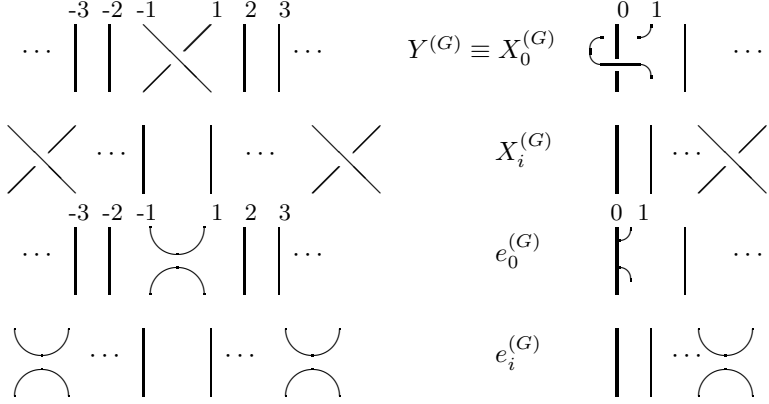
\begin{figure}[ht]
\begin{picture}(150,60)
\put(59,50){\mbox{$Y^{(G)}\equiv X_0^{(G)}$}}

\linethickness{0.2mm}
\put(2,50){\mbox{$\cdots$}}
\put(9,56){\mbox{{\small -3}}}
\put(14,56){\mbox{{\small -2}}}
\put(19,56){\mbox{{\small -1}}}
\put(40,56){\mbox{{\small 3}}}
\put(35,56){\mbox{{\small 2}}}
\put(30,56){\mbox{{\small 1}}}

\put(10,45){\line(0,1){10}}
\put(15,45){\line(0,1){10}}
\put(20,55){\line(1,-1){10}}
\put(20,45){\line(1,1){4}}
\put(26,51){\line(1,1){4}}
\put(35,45){\line(0,1){10}}
\put(40,45){\line(0,1){10}}
\put(42,50){\mbox{$\cdots$}}

\linethickness{0.4mm}
\put(90,45){\line(0,1){3}}
\put(90,50){\line(0,1){5}}
\linethickness{0.2mm}
\put(88,51){\oval(4,4)[l]}
\put(88,49){\line(1,0){5}}
\put(93,47){\oval(4,4)[tr]}
\put(93,55){\oval(4,4)[br]}
\put(95,56){\mbox{{\small 1}}}
\put(90,56){\mbox{{\small 0}}}

\put(100,45){\line(0,1){10}}

\put(107,50){\mbox{$\cdots$}}

\put(72,35){\mbox{$X_i^{(G)}$}}

\put(0,40){\line(1,-1){10}}
\put(0,30){\line(1,1){4}}
\put(6,36){\line(1,1){4}}
\put(45,40){\line(1,-1){10}}
\put(45,30){\line(1,1){4}}
\put(51,36){\line(1,1){4}}
\put(20,30){\line(0,1){10}}
\put(30,30){\line(0,1){10}}
\put(13,35){\mbox{$\cdots$}}
\put(35,35){\mbox{$\cdots$}}
\linethickness{0.4mm}
\put(90,30){\line(0,1){10}}
\linethickness{0.2mm}
\put(95,30){\line(0,1){10}}
\put(98,35){\mbox{$\cdots$}}
\put(102,40){\line(1,-1){10}}
\put(102,30){\line(1,1){4}}
\put(108,36){\line(1,1){4}}

\put(72,20){\mbox{$e_0^{(G)}$}}

\linethickness{0.2mm}
\put(2,20){\mbox{$\cdots$}}
\put(9,26){\mbox{{\small -3}}}
\put(14,26){\mbox{{\small -2}}}
\put(19,26){\mbox{{\small -1}}}
\put(40,26){\mbox{{\small 3}}}
\put(35,26){\mbox{{\small 2}}}
\put(30,26){\mbox{{\small 1}}}

\put(10,15){\line(0,1){10}}
\put(15,15){\line(0,1){10}}
\put(25,15){\oval(8,8)[t]}
\put(25,25){\oval(8,8)[b]}
\put(35,15){\line(0,1){10}}
\put(40,15){\line(0,1){10}}
\put(42,20){\mbox{$\cdots$}}

\linethickness{0.4mm}
\put(90,15){\line(0,1){10}}
\linethickness{0.2mm}
\put(90,17){\oval(4,4)[tr]}
\put(90,25){\oval(4,4)[br]}
\put(93,26){\mbox{{\small 1}}}
\put(89,26){\mbox{{\small 0}}}

\put(100,15){\line(0,1){10}}

\put(107,20){\mbox{$\cdots$}}

\put(72,5){\mbox{$e_i^{(G)}$}}
\linethickness{0.4mm}
\put(90,0){\line(0,1){10}}
\linethickness{0.2mm}
\put(95,0){\line(0,1){10}}
\put(98,5){\mbox{$\cdots$}}
\put(106,0){\oval(8,8)[t]}
\put(106,10){\oval(8,8)[b]}
\put(20,0){\line(0,1){10}}
\put(30,0){\line(0,1){10}}
\put(12,5){\mbox{$\cdots$}}
\put(32,5){\mbox{$\cdots$}}
\put(5,0){\oval(8,8)[t]}
\put(5,10){\oval(8,8)[b]}
\put(45,0){\oval(8,8)[t]}
\put(45,10){\oval(8,8)[b]}

\end{picture}
\caption{\label{generat} 
The graphical interpretation of the generators as symmetric tangles (on the left)
and as cylinder tangles (on the right)}

\end{figure}

$\ZB_n$ may be 
graphically interpreted (cf. figure \ref{generat}) as 
symmetric braids or cylinder braids: The symmetric picture
shows it as the group  of 
braids with $2n$ strands (numbered $-n,\ldots,-1,1,\ldots,n$) which are
fixed under a 180 degree rotation about the middle axis. 
In the cylinder picture one adds a single fixed line (indexed $0$)
on the left and
obtains $\ZB_n$ as the group of braids with $n$ strands that may 
surround this fixed line. 
The generators $X_i,i\geq0$ are mapped to the diagrams $X^{(G)}_i$ given in
figure \ref{generat}.

The group algebra of the B braid group   has  
 finite dimensional quotients which generalize the  
Temperley-Lieb, Hecke and Birman-Murakami-Wenzl algebras.
All these algebras support Markov traces that give rise to
B-type versions of the link polynomials of
Jones, HOMFLY and Kauffman.

  \unitlength0.9mm

 The Temperley-Lieb algebra of Coxeter type B 
has been introduced by tom Dieck
in \cite{tD1} as algebra of symmetric tangles without crossings.
\begin{de} The Temperley-Lieb Algebra ${\rm TB}_n$ of Coxeter type $B$ 
over a ring with
parameters $c,c',d$ is generated by $e_0,e_1,\ldots,e_{n-1}$ and relations
\begin{eqnarray}
e_1e_0e_1&=&c'e_1  \qquad e_0^2=de_0\qquad e_i^2=ce_i\\
e_ie_je_i&=&e_i\qquad |i-j|=1,i,j\geq 1\\
e_ie_j&=&e_je_i\qquad |i-j|>1
\end{eqnarray}
\end{de}
 The Hecke algebra $\HB_n$ has been studied by Dipper/James, tom Dieck,
Lambropoulou, Ariki and others.
\begin{de} \label{heckedef}
The Hecke Algebra of B-type  $\HB_n$ has generators
$X_0,X_1,\ldots,X_{n-1}$ over a ring with  parameters $Q,Q_0$ 
and relations (\ref{def1})-(\ref{def3}) and :
\begin{eqnarray}
X_0^2&=& (Q_0-1)X_0+Q_0
\qquad X_i^2= (Q-1)X_i+Q\quad i\geq 0\label{hdef4}
\end{eqnarray}
\end{de}
The reduced Birman-Wenzl algebra has been studied in \cite{rhobmw}.
\begin{de} Let $\BB_n$ 
be the algebra generated by invertible
$Y$, $X_1,\ldots,X_{n-1}$
over a ring with parameters
$\lambda$, $q$, $q_0:=q^{-1}$, $q_1$. 
The relations are: (\ref{def1})-(\ref{def2}) and:
\begin{eqnarray}
\delta&:=&q-q^{-1}\quad x:=1-\frac{\lambda-\lambda^{-1}}{\delta}
\quad
e_i:=1-\frac{X_i-X_i^{-1}}{q-q^{-1}}\\
X_ie_i&=&e_iX_i=\lambda e_i\qquad
e_iX_{i-1}^{\pm1}e_i=\lambda^{\mp1}e_i\\
X_1YX_1Y&=&YX_1YX_1\qquad YX_i=X_iY\qquad i>1\\
Y^2&=&q_1Y+q_0\qquad YX_1Ye_1=e_1
\end{eqnarray}
\end{de}
 \begin{satz} 
 \deu{Die Algebra} $\BB_n$\deu{ ist halbeinfach.}
\eng{is a semi simple algebra.}
\deu{Die einzelnen Komponenten
werden durch}\eng{The simple components are indexed by
the set  $\widehat{\Gamma}_n$
pairs of Young diagrams of size $n,n-2,\ldots$ }
\deu{ induziert}. \label{decomp}
\begin{equation}\BB_n=\bigoplus_{(\mu,\lambda)\in\widehat{\Gamma}_n}
 \BB_{n,(\mu,\lambda)}\end{equation}
 \deu{Die Bratteli-Regel lautet: Ein einfacher}
\eng{The Bratteli rule for restrictions of modules: A simple}
$\BB_{n,(\nu,\rho)}$ \deu{Modul}\eng{module} 
$V_{(\nu,\rho)},(\nu,\rho)\in\widehat{\Gamma}_n$ 
\deu{zerf\"allt in} \eng{decomposes into} $\BB_{n-1}$ 
\deu{Module zu}\eng{modules such that the $\BB_{n-1}$ module}  
$(\mu,\lambda)\in\widehat{\Gamma}_{n-1}$ 
\deu{genau dann, wenn}\eng{occurs iff}
$(\mu,\lambda)$ \deu{aus}\eng{may be obtained from} $(\nu,\rho)$ 
\deu{durch Wegnahme oder Erg\"anzung
einer Box entsteht}\eng{by adding or removing a box}.
 There exists a nondegenerate Markov trace on $\BB_n$.
\end{satz}
The dimension of $\BB_n$ is $2^n(2n-1)!!$. 	 In \cite{rhobmw}
a combinatorial basis was found.

  \unitlength1mm
 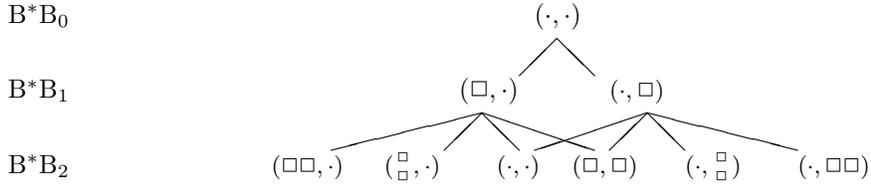
\begin{figure}[ht]
\begin{picture}(120,21)

\put(0,20){\mbox{$\BB_0$}}
\put(0,10){\mbox{$\BB_1$}}
\put(0,0){\mbox{$\BB_2$}}

\put(70,20){\mbox{$(\cdot,\cdot)$}}

\put(60,10){\mbox{$(\Box,\cdot)$}}
\put(80,10){\mbox{{\small $(\cdot,\Box)$}}}
\put(73,18){\line(1,-1){5}}
\put(73,18){\line(-1,-1){5}}

\put(65,0){\mbox{{\small $(\cdot,\cdot)$}}}
\put(75,0){\mbox{{\small $(\Box,\Box)$}}}
\put(35,0){\mbox{{\small $(\Box\Box,\cdot)$}}}
\put(50,0){\mbox{{\small $({\Box\atop\Box},\cdot)$}}}
\put(90,0){\mbox{{\small $(\cdot,{\Box\atop\Box})$}}}
\put(105,0){\mbox{{\small $(\cdot,\Box\Box)$}}}

\put(63,8){\line(-4,-1){20}}
\put(63,8){\line(1,-1){5}}
\put(63,8){\line(-1,-1){5}}
\put(63,8){\line(3,-1){15}}
\put(85,8){\line(-3,-1){15}}
\put(85,8){\line(1,-1){5}}
\put(85,8){\line(-1,-1){5}}
\put(85,8){\line(4,-1){20}}

\linethickness{0.2mm}

\end{picture}
\caption{\label{brat} 
\deu{Das Bratteli-Diagram von $\BB_n$}
\eng{The Bratteli digram of $\BB_n$}}

\end{figure}

Two dimensional integrable systems are described by solutions of the 
spectral parameter dependent Yang-Baxter-Equation (YBE) that reads:
$R_1(t_1)R_2(t_1t_2)R_1(t_2)=R_2(t_2)R_1(t_1t_2)R_2(t_1)$.
If the system is restricted to a half plane an additional matrix 
$K(t)\in{\rm End}(V)$  is needed to describe reflections. Is
has to fulfill Sklyanin's reflection equation:
\begin{equation}\label{re}
R(t_1/t_2)(K(t_1)\otimes1)R(t_1t_2)(K(t_2)\otimes1)=
(K(t_2)\otimes1)R(t_1t_2)(K(t_1)\otimes1)R(t_1/t_2)
\end{equation}
Solutions of the Yang Baxter equation can be
 obtained from the standard (type A)
Birman-Wenzl algebra by the following Baxterization 
procedure:
\begin{equation}
R_i(t)= -\delta t (t+q\lambda^{-1})+
(t-1)(t+q\lambda^{-1})X_i+\delta t(t-1)e_i
\end{equation}
\begin{satz} $K(t)=(t^2q_1(1-t^2)^{-1}+Y)f_1(t)$ is  (for all $f_1$) a solution of
the reflection equation (\ref{re}).
\end{satz}

The proof of this result is given in \cite{rhoref}.

$\BB_n$ supports a Markov trace that can be used to 
define a link invariant for 
links of B-type which are links in a solid torus. 
There is an analog of Markov's theorem for type B links found by S.
Lambrodopoulou in \cite{lamb}. It 
 implies that there exists an extension of the Kauffman polynomial 
to braids of B-type. For a B-type link $\hat{\beta}$
that is the closure of a B-braid $\beta_\in\ZB_n$ we define:
\begin{de}
The B-type Kauffman polynomial of a B-link $\hat{\beta}$ is defined to be
\begin{equation}
L(\hat{\beta},n):=x^{n-1}\lambda^{e(\beta)}\tr(\beta)\qquad\beta\in\ZB_n
\end{equation}
$e:\ZB_n\rightarrow\ZZ$ is the exponential sum with $e(X_i)=1,e(Y)=0$.
\end{de}

\section{The Potts-Model with a boundary}

\begin{figure}[ht]
\begin{picture}(110,36)
\linethickness{0.5mm}

\put(20,10){\circle*{2}}
\put(20,20){\circle*{2}}
\put(20,30){\circle*{2}}
\put(30,10){\circle*{2}}
\put(30,20){\circle*{2}}
\put(30,30){\circle*{2}}
\put(40,10){\circle*{2}}
\put(40,20){\circle*{2}}
\put(40,30){\circle*{2}}

\put(20,10){\line(0,1){20}}
\put(30,10){\line(0,1){20}}
\put(40,10){\line(0,1){20}}

\put(20,10){\line(1,0){20}}
\dottedline{10}{10}{1}{0}{10}
\put(20,20){\line(1,0){20}}
\dottedline{10}{20}{1}{0}{10}
\put(20,30){\line(1,0){20}}
\dottedline{10}{30}{1}{0}{10}

\multiput(5,4)(0,3){10}{\line(1,1){5}}

\linethickness{1mm}
\put(10,5){\line(0,1){30}}

\linethickness{0.2mm}

\put(20,35){\line(1,-1){8}}
\put(21,24){\line(1,-1){8}}
\put(21,14){\line(1,-1){9}}
\put(30,35){\line(1,-1){8}}
\put(31,24){\line(1,-1){8}}
\put(31,14){\line(1,-1){9}}

\put(16,11){\line(1,1){8}}
\put(16,21){\line(1,1){8}}
\put(26,11){\line(1,1){8}}
\put(26,21){\line(1,1){8}}
\put(36,11){\line(1,1){9}}
\put(36,21){\line(1,1){9}}

\put(20,5){\line(1,1){4}}
\put(30,5){\line(1,1){4}}

\put(26,31){\line(1,1){4}}
\put(36,31){\line(1,1){4}}

\put(41,14){\line(1,-1){4}}
\put(41,24){\line(1,-1){4}}

\put(40,35){\line(1,-1){5}}

\put(40,5){\line(1,1){5}}

\linethickness{0.4mm}
\put(15,0){\line(0,1){7}}
\put(15,9){\line(0,1){8}}
\put(15,19){\line(0,1){8}}
\put(15,29){\line(0,1){8}}
\linethickness{0.2mm}
\put(14,9.5){\oval(3,3)[l]}
\put(16,8){\line(1,-1){3.5}}
\put(14,19.5){\oval(3,3)[l]}
\put(16,18){\line(1,-1){3.5}}
\put(14,29.5){\oval(3,3)[l]}
\put(16,28){\line(1,-1){3.5}}
\put(20,35){\line(-1,-1){4}}

\end{picture}
\caption{\label{latticelink} A lattice with boundary 
together with its graph}
\end{figure}
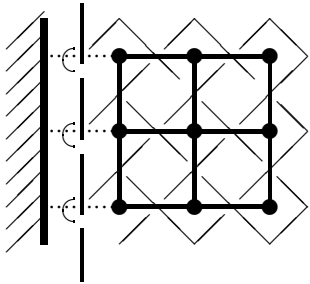

We investigate  a generalization of the ordinary
Potts  model \cite{kauff} by including a  reflecting boundary.
Hence we have besides the usual lattice of sites (denoted by $V$) a wall
interacting with it.
An example is shown in figure \ref{latticelink}. The dotted lines
indicate interaction bonds with the wall while the solid lines are
usual bonds between sites. Each site supports a 'spin' which
may occupy one of $f$ states.

The set of all states is 
${\cal S}=\{S:V\rightarrow\{0,1,\ldots,f-1\}\}$.
The partition function (with $k$ being the Boltzmann constant and $T$
the temperature)   is
\begin{equation}\label{zdef}
Z_G=\sum_{S\in{\cal S}}{\rm exp}\left(\frac{-E(S)}{kT}\right)
\end{equation}
with the Hamiltonian  ($\delta(x,y)$ is the Kronecker symbol
with values $0,1$)
\begin{equation}
E(S)=\sum_{(i,j)\in B_1}\delta(S_i,S_j)+
\kappa\sum_{i\in B_0}(1-\delta(0,S_i))
\end{equation}

Here $B_0\subset V$ is the set of sites which have boundary bonds, and
$B_1$ is the set of inner bonds. 
The first term in $E(S)$ is the usual Hamiltonian of the Potts model.
The second term introduces the boundary condition.

In analogy with Kauffman's treatment of the ordinary Potts
 model \cite{kauff}
we associate a link diagram  of type B with the boundary lattice. 
The partition function may then be expressed as a normalization of 
the Markov trace of this link.

\section{Coxeter-B braided categories}

The language of braided tensor categories (BTCs) turned out to be 
the right framework for braid group applications from knot theory to 
QFT. In this section we generalize it to Coxeter type B. 
Details are given in \cite{rhocat}.

\begin{de} Let $\kat$ be a  rigid ribbon BTC. 
A  Coxeter-B braided category over $\kat$ is 
an embedding of $\kat$ in a  rigid monoidal category $\widehat{\kat}$
which has the same objects
and obeys the following list of axioms hold.
Morphisms $\Mor(X,Y):=\Mor_\kat(X,Y)$ are said to be local and 
morphisms $\Mor^{(G)}(X,Y):=\Mor_{\widehat{\kat}}(X,Y)$ are said 
to be global.
\begin{eqnarray}
\forall X&&\exists b_X\in\Mor^{(G)}(X,X)\\
b_Yf&=&fb_X\qquad\forall f\in\Mor(X,Y)\label{bycommu}\\
\id_X\otimes b_Y&=&\Psi_{Y,X}(b_Y\otimes\id_X)\Psi_{X,Y}\label{y2form}\\
b_X\otimes b_Y&=&b_{X\otimes Y}\Psi_{Y,X}\Psi_{X,Y}
  =\Psi_{Y,X}\Psi_{X,Y}b_{X\otimes Y}\label{bytp}\\
b_X^\ast&=&\sigma(X^\ast)^2b_{X^\ast}^{-1} \label{bydual}
\end{eqnarray}
\end{de}
Note that $\Psi$ is a braiding of $\kat$, not of $\widehat{\kat}$.
This makes (\ref{y2form}) possible which  otherwise would give a
contradiction to naturality of $\Psi$. The axioms imply: 
\begin{eqnarray} 
 \Psi_{Y,X}(b_Y\otimes\id_X)\Psi_{X,Y}(b_X\otimes\id_Y)&=&
 (b_X\otimes\id_Y)\Psi_{Y,X}(b_Y\otimes\id_X)\Psi_{X,Y}	 \label{lem1b} \\
\ev_X(b_{X^\ast}\otimes\id_X)\Psi_{X,X^\ast}(b_X\otimes\id_{X^\ast})&=&
\ev_X\Psi_{X,X^\ast}\label{lem1c}
\end{eqnarray}
A restricted Coxeter-B braided tensor category is defined similarly, but
with 
(\ref{y2form}) replaced by the two relations above.
The graphical origin of this axioms becomes clear by drawing
pictures. A nice example of a Coxeter-B braided category is given by
the category of Aplimorphisms of a quasitriangular Hopf algebra.

\begin{de} A  restricted Coxeter-B braided ((weak) quasi) Hopf 
algebra $H$ is a
 ((weak) quasi) quasi triangular ribbon Hopf algebra 
 with an element $\overline{v}\in H$
 such that
 \begin{eqnarray}
    R_{2,1}\overline{v}_2R\overline{v}_1&=&
 \overline{v}_1R_{2,1}\overline{v}_2R	
\qquad   \epsilon(\overline{v})=1\label{reshd1}\label{resvepsilon}\\
 \Delta(\overline{v})&=&R^{-1}(1\otimes\overline{v})
 R(\overline{v}\otimes1)\label{reshd2}
 \end{eqnarray}
\end{de}

Such algebras can be constructed explicitly from the quantum
Weyl group \cite{DH}.

We have the following Tannaka-Krein style duality between 
Coxeter B-type braided Hopf algebras and B-type tensor categories.

\begin{satz} The representation category $\Rep(H)$ of 
a restricted Coxeter-B braided ((weak) quasi) Hopf algebra algebra is a 
restricted B-braided tensor category.

Contrary, if $\widehat{\kat}$ is a
Coxeter-B braided tensor category 
and $F:\widehat{\kat}\rightarrow\Vec$ is a 
((weak) quasi) tensor functor in the restricted sense that
the naturality equation $c_{X,Y}(F(f)\otimes F(g))=F(f\otimes g)c_{X,Y}$
holds only if $g$ is a local morphism then 
the set ${\rm Nat}(F,F)$ of endo transformations that are natural
for local morphism is a 
restricted Coxeter-B braided  ((weak) quasi) Hopf algebra.
\end{satz}

Consider a quantum field theory \cite{haag} specified by a net of local
observables $\oalg(\kegel)$ living on the half plane
$\{(x_1,x_2)\in\RR^2\mid x_1\geq 0\}$. 
We assume that boundary conditions are imposed in such a way that we have 
reflection at the line $(0,\RR)$ by letting the full translation group
$\RR^2$ act on the half plane. This action is not free and this will lead to 
global morphisms and hence to the occurrence of a Coxeter-B braided tensor 
category.

\begin{figure}[ht]
\begin{picture}(150,60)

\linethickness{0.4mm}
\put(5,0){\line(0,1){60}}
\linethickness{0.1mm}
\multiput(0,0)(0,3){19}{\line(1,1){5}}
\linethickness{0.2mm}

\put(15,45){\line(1,1){10}}
\put(15,45){\line(1,-1){10}}
\put(35,45){\line(-1,1){10}}
\put(35,45){\line(-1,-1){10}}
\put(23,45){\mbox{${\mathcal O}_1$}}
\put(15,45){\line(-1,1){10}}
\put(15,45){\line(-1,-1){10}}
\put(6,45){\mbox{${\mathcal O}_1'$}}
\put(39,45){\mbox{${\mathcal O}_1'$}}
\put(35,45){\line(1,1){15}}
\put(35,45){\line(1,-1){15}}

\put(17,25){\line(-1,1){8}}
\put(17,25){\line(-1,-1){8}}
\put(8,23){\mbox{${\mathcal O}_2$}}
\put(9,33){\line(-1,-1){4}}
\put(9,17){\line(-1,1){4}}

\put(5,0){\line(1,1){8}} \put(5,16){\line(1,-1){8}}
\put(6,7){\mbox{${\mathcal O}_3$}}

\linethickness{0.4mm}
\put(80,0){\line(0,1){60}}
\linethickness{0.1mm}
\multiput(75,0)(0,3){19}{\line(1,1){5}}
\linethickness{0.2mm}

\put(85,10){\line(1,1){10}}
\put(85,10){\line(1,-1){10}}
\put(105,10){\line(-1,1){10}}
\put(105,10){\line(-1,-1){10}}
\put(93,10){\mbox{${\mathcal O}_1$}}

\put(125,30){\line(1,1){10}}
\put(125,30){\line(1,-1){10}}
\put(145,30){\line(-1,1){10}}
\put(145,30){\line(-1,-1){10}}
\put(133,30){\mbox{${\mathcal O}_2$}}

\put(110,10){\vector(2,1){20}}
\put(120,10){\mbox{$U$}}

\put(90,18){\vector(-4,1){10}}
\put(80,22){\vector(4,1){40}}
\put(100,22){\mbox{$V$}}

\end{picture}
\caption{The half plane with various double cones (on the left) and with
the reflecting transportation that leads to global intertwiners (on the right)
\label{halfplane}} 
\end{figure}
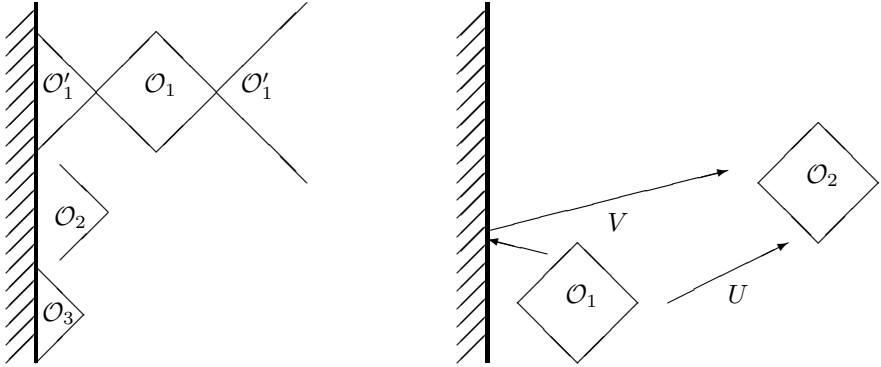

Fields shall be localized in double cones. Here we extend the usual notion
of a double cone to include all translations of double cones. Thus we also 
have regions like those in the left of figure \ref{halfplane}. 
This figure also shows the causal complement $\kegel'$ of a double cone.
A double cone that does not touch the boundary shall be called 
regular.
Further we assume isotony and  locality and the existence of a vacuum
representation $\pi_0$ which is translation invariant and faithful for
all local algebras $\oalg(\kegel)$ of regular double cones $\kegel$.

Now, let $\kegel_1,\kegel_2$ be two causally disjoint double cones 
of equal size as
shown in the right half of figure \ref{halfplane}. 
Further let $\rho_1$ be a transportable 
morphisms localized in $\kegel_1$ and  let $\rho_2\sim\rho_1$ be
localized in $\kegel_2$. Assume that $\rho_1$ (and thus $\rho_2$)
is irreducible in the sense that $\pi_0\circ\rho_1$ is an irreducible 
representation of $\oalg$.
There are two translations that map $\kegel_1$ onto $\kegel_2$:
A direct one and one that passes the reflecting boundary. 
Thus we have two charge transporters $U,V\in\Mor^{(G)}(\rho_1,\rho_2):=
\{T\in\oalg\mid T\rho_1(A)=\rho_2(A)T\forall A\in\oalg\}$.
Hence we have a self intertwiner $Y:=U^{-1} V\in\Mor^{(G)}(\rho_1,\rho_1)$.
We see that the vacuum representation may not be faithful in the presence of a 
boundary. The localized and transportable morphisms form a 
Coxeter-B braided tensor category.

  \unitlength0.8mm

\end{document}